\documentclass{article}
\usepackage{amssymb,latexsym,amsmath,amsthm,mathrsfs}
\usepackage{graphicx}
\newcommand{\comment}[1]{}
\begin{document}
\title{On highly transcendental quantities which cannot be expressed by
integral formulas\footnote{Presented to the St. Petersburg Academy on
October 16, 1775.
Originally published as
{\em De plurimis quantitatibus transcendentibus quas nullo modo
per formulas integrales exprimere licet},
Acta academiae scientiarum Petropolitanae \textbf{4} (1784), no. II,
31--37.
E565 in the Enestr{\"o}m index.
Translated from the Latin by Jordan Bell,
Department of Mathematics, University of Toronto, Toronto, Canada.
Email: jordan.bell@utoronto.ca}}
\author{Leonhard Euler}
\date{}
\maketitle

1. Integral formulas, whose integration cannot be obtained in terms of
algebraic quantities, are commonly regarded as the single source of
all transcendental quantities. Thus from the integral formula
\[
\int \frac{dx}{x}
\]
are born the logarithms, and from the formula
\[
\int \frac{dx}{1+xx}
\]
circular arcs; these quantities, although transcendental, have now indeed been
absorbed into Analysis, so that they can be dealt with as easily as algebraic
quantities. Besides these indeed are also the quantities which involve
the rectification of conic sections, which now have been so explored by Geometers that problems which lead to these can obtain perfect solutions. These
transcendental quantities are further contained in the integral formulas
of the type\footnote{Translator: This is an elliptic integral.
Euler's 1760 paper ``Consideratio formularum,
quarum integratio per arcus sectionum conicarum absolvi potest'', E273, looks
like it's about integrals
of this form.}
\[
\int dx \surd \frac{f+gxx}{b+kxx}.
\]
Whenever other transcendental quantities occur, they always permit
representation as the quadrature of certain curved lines; doing
this leads to integral formulas which, though complicated, are able to express
the true values of such transcendental quantities.

2. I have observed nevertheless that innumerable types of transcendental quantities
can be exhibited which cannot be expressed in any way by integral formulas,
even if their true values can often be defined easily enough. Quantities
of this type arise principally from infinite series  
whose sums can so far be reduced in no way to integral formulas.
Especially noteworthy among these is this infinite series, differing slightly
from the geometric,
\[
1+\frac{1}{3}+\frac{1}{7}+\frac{1}{15}+\frac{1}{31}+\frac{1}{63}+\textrm{etc.},
\]
the denominators of whose fractions are powers of two less unity, whose
approximate sum is indeed not very difficult to assign. It is easily
understood moreover that this can be expressed neither by a rational nor by irrational numbers; then indeed it also seems certain enough that it can be
expressed neither by logarithms or circular arcs. 
So far however no path is apparent for investigating an integral formula of this
type, whose value exhibits exactly the sum of this series.
But letting all the terms of this series each be converted into a geometric
series in the usual way,
it is worth noting that its sum can be represented by the following formulas
\[
\int \frac{1}{2^{xx}}+2\int\frac{1}{2^{xx+x}}+2\int\frac{1}{2^{xx+2x}}
+2\int\frac{1}{2^{xx+3x}}+\textrm{etc.},
\]
where $\int \frac{1}{2^{xx}}$ denotes the sum of the infinite series whose
term corresponding to the index $x$ is $\frac{1}{2^{xx}}$, which will therefore
be
\[
\frac{1}{2}+\frac{1}{2^4}+\frac{1}{2^9}+\frac{1}{2^{16}}+\frac{1}{2^{25}}+
\textrm{etc.};
\]
but equally as the sum of this series cannot be exhibited in any way, the same
should be understood for the remaining parts. The true sum of the given series
is approximately $1,606695152$; if this number were found to be equal to
any known number, such as
a logarithm or a circular arc, 
it would without doubt be an extraordinary discovery.

3. Like how we have subtracted unity from the powers of two here, thus
let us add these to the indefinite number $x$, and let us put
\[
y=\frac{1}{1+x}+\frac{1}{2+x}+\frac{1}{4+x}+\frac{1}{8+x}+\textrm{etc.},
\]
which equation, if $x$ is taken for the abscissa and $y$ for the ordinate,
expresses a certain curved line in which in fact one can assign a single
ordinate corresponding to the abscissa $x=0$, which of course will be $y=2$.
But for all the other abscissas, the ordinates will be highly
transcendental quantities, which thus far seem not to be able to be expressed
by integral formulas, so that the nature of this curve is such that it cannot
be expressed by any equation, either differential or integral. In the meanwhile
however it is clear that the abscissas
\[
x=-1, \quad x=-2, \quad x=-4, \quad x=-8, \quad x=-16 \quad \textrm{etc.}
\]
lead to infinitely large ordinates, while on the other hand taking $x=\infty$
the ordinates will vanish.

4. This equation can be made more general if we let any particular number $a$
be assumed in place of $2$, so that it would be
\[
y=\frac{1}{1+x}+\frac{1}{a+x}+\frac{1}{a^2+x}+\frac{1}{a^3+x}+\frac{1}{a^4+x}+\textrm{etc.},
\]
where the ordinate corresponding to the abscissa $x=0$ is
\[
y=\frac{a}{a-1},
\]
and indeed again all the remaining terms are highly transcendental. It is
also clear that if one takes $a=1$, all the ordinates will then be infinite,
unless perhaps the abscissa $x$ is taken as infinite. Like how we gave
all the fractions here the $+$ sign, they can also thus alternate, so that
it would be
\[
y=\frac{1}{1+x}-\frac{1}{a+x}+\frac{1}{a^2+x}-\frac{1}{a^3+x}+\frac{1}{a^4+x}-\textrm{etc.};
\]
while now the ordinate that will correspond to the abscissa $x=0$ is
\[
y=\frac{a}{a+1}.
\]
Moreover, in place of $x$ its powers could be written in the following terms,
so that an equation of this type would be obtained
\[
y=\frac{1}{1+1}+\frac{1}{a+x}+\frac{1}{a^2+x^2}+\frac{1}{a^3+x^3}+\frac{1}{a^4+x^4}+\textrm{etc.},
\]
which certainly is composed such that the nature of this curve seems not to be
able to be
expressed by any finite equation, either differential or integral.

5. Truly besides these forms there are infinitely many others which
can be exhibited, proceeding
only by powers of some quantity $x$, which also cannot be reduced to integral
formulas in any way. 
For since so far no series of this type have been able to be summed unless their
exponents progress in a geometric progression,
as soon as we take any other law for the progression of the exponents of $x$
their summation always seems to overpower the strength of analysis,
even if we call on the aid of complicated integral formulas.
Thus if the  infinite equation
\[
y=1+x+x^3+x^6+x^{10}+x^{15}+x^{21}+\textrm{etc.},
\]
is given, where the exponents of $x$ are the trigonal numbers,
altogether no finite equation can be exhibited for such a curve by just integral formulas.
Here indeed it at once clear that if it were $x=1$ or $>1$, the ordinates
will increase to infinity; truly if we take for $x$ values less than unity,
the true values of the ordinate $y$
can be assigned approximately easily enough; for instance if we take $x=\frac{1}{10}$ we will have at once in decimal fractions
\[
y=1,1010010001000010000010000001,
\]
which fraction can be continued as far as we desire very easily, simply
by continually increasing by unity the number of zeros placed between the unities.
While the true value of this fraction cannot be exhibited by integral formulas,
we can still form an idea exact enough for us, so that if
perhaps the quadrature of the circle could be reduced to such a form,
it would be found as thoroughly as could be obtained;
and the same will hold for all the other decimal fractions all of whose
figures progress according to a certain law, which however contain no repeating
periods, in which cases of course the value could be exhibited as a rational.

6. It would without doubt be most momentous if the sum of this series
\[
1+x+x^3+x^6+x^{10}+x^{15}+\textrm{etc.}
\]
could be exhibited in a finite expression, however transcendental. 
For from this a very coherent demonstration of the theorem of Fermat could be approached,
that all integral numbers are the sum of three trigonal numbers;
for the problem would just be to resolve the cube of this sum into a series
and then to show that all the powers of $x$ must occur. 

7. Since here we took $x=\frac{1}{10}$, so the expansion into a decimal
fraction works out easily, the approximate summation
of that series
will indeed be hardly more difficult
if any fraction less than unity is taken for $x$. For instance if we take
$x=\frac{1}{2}$, so that it becomes
\[
y=1+\frac{1}{2}+\frac{1}{2^3}+\frac{1}{2^6}+\frac{1}{2^{10}}+\frac{1}{2^{15}}
+\frac{1}{2^{21}}+\textrm{etc.},
\]
by a decimal fraction it will be approximately
\[
y=1,64163256066.
\]

8. Moreover all series of this type can be included under this general form
\[
y=Ax^\alpha+Bx^\beta+Cx^\gamma+Dx^\delta+\textrm{etc.},
\]
providing the exponents $\alpha,\beta,\gamma,\delta$ etc. do not progress in
arithmetic ratio. For whatever law is set between the coefficients $A,B,C,D$ etc.,
it is certain that at no time can the value of $y$ be expressed in a finite
way,
no matter what quadratures are called upon for aid. 
And in turn if forms of a similar nature are multiplied, divided or otherwise
combined in any way, types of transcendental quantities of this kind
will always arise, which seem to elude all integral formulas. Neither
indeed do these new types of transcendentals admit comparison between each other,
whence we cannot admire enough the immense multitude of all the different types
of quantities.

9. For since between two numbers, no matter how close, not only innumerable
rational fractions but also infinitely many irrationals can be assigned, all
of which
are different not only from the former but to each other,
it seem scarcely believable for there also to be given 
infinitely many other types of quantities which differ among themselves and between
themselves.
And since now such an infinity of types is demonstrated, which are
able to be exhibited by integral formulas,
we are compelled to recognize, even besides the types given so far,
innumerably many types of transcendental quantities which do not themselves
allow reduction to the previous in any way.

\end{document}